\documentclass[a4,11pt]{article}%
\usepackage{amsmath}
\usepackage{amsfonts}
\usepackage{amssymb}
\usepackage{graphicx}%
\providecommand{\U}[1]{\protect\rule{.1in}{.1in}}
\providecommand{\U}[1]{\protect\rule{.1in}{.1in}}
\providecommand{\U}[1]{\protect\rule{.1in}{.1in}}
\RequirePackage[colorlinks,citecolor=blue,urlcolor=blue]{hyperref}
\RequirePackage{hypernat}
\providecommand{\U}[1]{\protect\rule{.1in}{.1in}}
\numberwithin{equation}{section}
\newtheorem{lemma}{Lemma}[section]
\newtheorem{theorem}{Theorem}[section]
\newtheorem{corollary}{Corollary}[section]
\newtheorem{proposition}{Proposition}[section]
\newtheorem{definition}{Definition}[section]

\newtheorem{remark}{Remark}[section]
\begin{document}

\title{\vspace{-1in}\parbox{\linewidth}{\footnotesize\noindent
} \vspace{\bigskipamount} \\A complex Feynman-Kac formula via linear backward stochastic differential
equations
\thanks{The work is supported by Natural Science Foundation of China and Jiangsu
province(No.11401414; No.BK20140299; No.14KJB110022) and the collaborative
innovation center for quantitative calculation and control of financial
risk. } }
\author{Yuhong Xu\\Mathematical center for interdiscipline research and\\School of mathematical sciences, Soochow University,\\Suzhou 215006, P. R. China.\\E-mail: yuhong.xu@hotmail.com}
\date{}
\maketitle

\indent\textbf{Abstract.} A complex notion of backward stochastic
differential equation (BSDE) is proposed in this paper to give a
probabilistic interpretation for linear first order complex partial
differential equation (PDE). By the uniqueness and existence of regular
solutions to complex BSDE, we deduce that there exists a unique classical
solution $\{\mathbb{U}(t,x)$ to complex PDE and $\{\mathbb{U}(t,x)$ is
analytic in $x$ for each $t$. Thus we extend the well known real Feynman-Kac
formula to a complex version. It is stressed that our complex BSDE
corresponds to a linear PDE without the second order term.


\medskip


\section{Introduction}

The Feynman-Kac formula, named after Richard Feynman and Mark Kac,
establishes a link between PDEs and stochastic differential equations
(SDEs). It offers a method of solving certain PDEs by simulating random
paths of a stochastic process. Linear real-valued Feynman-Kac formula was
studied early in Kac (1949, 1951) as a formula for determining the
distribution of certain Wiener functionals. Then by the theory of BSDE,
Pardoux and Peng (1992) and Peng (1991, 1992) generalized them to nonlinear
versions. They also derived some stochastic versions (Pardoux and Peng 1992,
Peng 1992). There are many other papers in this direction, however we don't
list them all here. The present paper propose a complex notion of BSDE and
deduce a complex Feynman-Kac formula for linear first order complex PDE.

Let us be more precise. We first introduce a one dimensional complex BSDE:

\begin{equation}
\mathbb{Y}_{s}^{t,x}=h(\mathbb{X}_{T}^{t,x})+\int_{s}^{T}g\left( r,\mathbb{X}%
_{r}^{t,x},\mathbb{Y}_{r}^{t,x},\mathbb{Z}_{r}^{t,x},\mathbb{T}%
_{r}^{t,x}\right) dr-\int_{s}^{T}\mathbb{Z}_{r}^{t,x}d\mathbb{B}%
_{r}-\int_{s}^{T}{}\mathbb{T}_{r}^{t,x}d\bar{\mathbb{B}}_{r},\ t\leq s\leq T.
\label{1.1}
\end{equation}%
where ${(}\mathbb{B}{{_{s})}_{s\in \lbrack 0,T]}}$ is a complex Brownian
motion with ${(}\bar{\mathbb{B}}{{_{s})}_{s\in \lbrack 0,T]}}$ its conjugate
counterpart, $\left( \mathbb{X}^{t,x}\right) $ is defined as%
\begin{equation}
\mathbb{X}_{s}^{t,x}=x+\int_{t}^{s}\sigma \left( r\right) d\mathbb{B}%
_{r}+\int_{t}^{s}\gamma \left( r\right) d\bar{\mathbb{B}}_{r},\ t\leq s\leq
T.  \label{1.2}
\end{equation}%
We next want to find a triple of adapted processes $\{(\mathbb{Y}_{s}^{t,x},%
\mathbb{Z}_{s}^{t,x},\mathbb{T}_{s}^{t,x});t\leq s\leq T\}$ with values in $%
\mathbb{C\times C\times C}$ which solves uniquely \eqref{1.1}. We finally
show that under some analytic conditions on coefficients, $\{\mathbb{U}(t,x)=%
\mathbb{Y}_{t}^{t,x};0\leq t\leq T,x\in \mathbb{C}\}$ is the unique solution
of the following complex PDE:

\begin{equation}
\left\{ \begin{aligned} \Bbb{U}_t(t,x) &= -g\left( t,x,\Bbb{U}(t,x),\sigma
_t \Bbb{U}_x(t,x),\gamma _t \Bbb{U}_x(t,x)\right), \\ \Bbb{U}(T,x)&=h(x),\
0\leq t\leq T. \end{aligned}\right.
\end{equation}%
where $\mathbb{U}$ takes values in $\mathbb{C}$ and is analytic with respect
to $x$ for each $t$. Note that not like the real Feynman-Kac formula , BSDE %
\eqref{1.1} gives a probabilistic interpretation to PDEs without the second
order term, which is essentially due to the complex It\^{o}'s formula
involving analytic functions. We refer to Ub$\emptyset $e (1987), Davis
(1979) and Varopoulos (1981) for complex stochastic analysis.

This paper is organized as follows. In section 2, we make some
preliminaries. Section 3 proves existence, uniqueness and regularity for the
solutions of complex BSDEs. In section 4, we establish a link between a
class of linear PDEs and complex BSDEs.

\section{Preliminaries}

Let ${(\mathbf{B}_{t})}=(B_{t}^{1},B_{t}^{2}){_{t\geq 0}}$ be a standard
Brownian motion in $\mathbf{R}^{2}$ on a probability space $(\Omega ,%
\mathcal{F},\mathit{P})$ and $(\mathbb{B}_{t}){_{t\geq 0}}$ be its complex
counterpart, i.e. $\mathbb{B}_{t}=B_{t}^{1}+iB_{t}^{2}$ where $i=\sqrt{-1}$
is the imaginary unit. Let $(\mathcal{F}_{s}^{t}){_{s\geq t}}$ be the
augmented Brownian filtration generated by the Brownian motion ${(\mathbf{B}%
_{s})}_{s\geq t}$ from time $t$. $T<0$ is a fixed time. Throughout the paper
we will work within the time interval $[0,T]$. For $x\in \mathbf{R}^{2}$, $x$
means a column vector $\left(
\begin{array}{l}
x^{1} \\
x^{2}%
\end{array}%
\right) .$ We write $x=x^{1}+ix^{2}$ as its complex counterpart in $\mathbb{C%
}$, the set of complex numbers. For $x\in \mathbf{R}^{2}$ or $x\in \mathbb{C}
$, We define a common Euclid norm $|{x}|=\sqrt{\ |x^{1}|^{2}+|x^{2}|^{2}}$.
For $x,y\in \mathbf{R}^{2}$, we denote by $\left\langle x,y\right\rangle $
the scalar product of $x,y$. Let $\mathbf{M,N}$ be two fields of (real or
complex) numbers. We denote by ${\mathcal{H}_{\mathcal{F}}^{\mathrm{2}}}(t,T;%
\mathbf{N})$ the space of all $\mathcal{F}_{s}^{t}$-progressively measurable
$\mathbf{N}$-valued processes $\{(\varphi _{s});t\leq s\leq T\}$ s.t. $%
\mathbf{E}\left[ \int_{t}^{T}|\varphi _{s}|^{2}ds\right] <\infty $ and by ${%
\mathcal{S}}_{\mathcal{F}}^{\mathrm{2}}\left( t,{\mathit{T}};\mathbf{N}%
\right) $ the set of {continuous and progressively measurable }$\mathbf{N}$%
-valued processes $\{(\psi _{s});t\leq s\leq T\}$ s.t. ${\mathbf{E}}\left[ {%
\mathrm{sup}}_{t\leq s\leq T}|\psi _{s}|^{2}\right] <\infty \}.$ $C^{k}(%
\mathbf{M,N}),C_{b}^{k}(\mathbf{M,N}),C_{p}^{k}(\mathbf{M,N})$ will denote
respectively the set of functions of class $C^{k}$ from $\mathbf{M}$\textbf{%
\ }to\textbf{\ }$\mathbf{N}$ , the set of those functions of class $C^{k}$
whose partial derivatives of order less than or equal to $k$ are bounded,
and the set of those functions of class $C^{k}$ which, together with all
their partial derivatives of order less than or equal to $k$, grow at most
like a polynomial function of the variable $x$ at infinity.

For a single-variable function $f$, $f^{\prime }$ denotes its derivative,
for a multi-variable function $f$, we sometimes denote by $f_{x}$ its
partial derivative w.r.t $x$ variable.

We now announce a result about $2\times 2$ matrices.

\begin{definition}
A real $2\times 2$ matrix $A$ is of class $\mathbb{C}_L$ if and only if it
has the form
\begin{equation}
\left(
\begin{array}{ll}
a & -b \\
b & a%
\end{array}
\right)
\end{equation}

where $a,b\in \mathbf{R}.$

By a direct calculation, we have
\end{definition}

\begin{lemma}
\label{lem2.2}Class $\mathbb{C}_{L}$ is an exchangeable semigroup, i.e.: let
$2\times 2$ matrices $A$, $B$ be of class $\mathbb{C}_{L}$, then{\normalsize %
\ }$A+B${\normalsize , }$A\cdot B${\normalsize = }$B\cdot A${\normalsize , }$%
\lambda A${\normalsize , (}$\lambda \in \mathbf{R}${\normalsize )} are all
of class $\mathbb{C}_{L}$. If $A^{-1}$ exists, then $A^{-1}$ also belongs to
class $\mathbb{C}_{L}$.

Proof. The following calculation leads to the above results:

$\left(
\begin{array}{ll}
a & -b \\
b & a%
\end{array}
\right) +\left(
\begin{array}{ll}
c & -d \\
d & c%
\end{array}
\right) =\left(
\begin{array}{ll}
a+c & -b-d \\
b+d & a+c%
\end{array}
\right) $;

$\left(
\begin{array}{ll}
a & -b \\
b & a%
\end{array}
\right) \left(
\begin{array}{ll}
c & -d \\
d & c%
\end{array}
\right) =\left(
\begin{array}{ll}
ac-bd & -ad-bc \\
ad+bc & ac-bd%
\end{array}
\right) =\left(
\begin{array}{ll}
c & -d \\
d & c%
\end{array}
\right) \left(
\begin{array}{ll}
a & -b \\
b & a%
\end{array}
\right) $;

$\left(
\begin{array}{ll}
a & -b \\
b & a%
\end{array}%
\right) ^{-1}=\frac{1}{a^{2}+b^{2}}\left(
\begin{array}{ll}
a & b \\
-b & a%
\end{array}%
\right) $ .
\end{lemma}

\section{Complex BSDE: existence, uniqueness and regularity}

Let $\mathbb{X}^{t,x}$ satisfy \eqref{1.2} with $\sigma \left( \cdot \right)
,\gamma \left( \cdot \right) $ being deterministic and square-integrable
functions. We introduce the following complex BSDE:

\begin{equation}
\mathbb{Y}_{s}^{t,x}=h(\mathbb{X}_{T}^{t,x})+\int_{s}^{T}g\left( r,\mathbb{X}%
_{r}^{t,x},\mathbb{Y}_{r}^{t,x},\mathbb{Z}_{r}^{t,x},\mathbb{T}%
_{r}^{t,x}\right) dr-\int_{s}^{T}\mathbb{Z}_{r}^{t,x}d\mathbb{B}%
_{r}-\int_{s}^{T}{}\mathbb{T}_{r}^{t,x}d\bar{\mathbb{B}}_{r},\ t\leq s\leq T
\label{3.1}
\end{equation}%
where $\mathit{h}$ : $\mathbb{C\rightarrow C}$, $\mathit{g=}\widetilde{g}%
\left( r,x\right) +\alpha \left( r\right) y+\beta \left( r\right) z+\theta
\left( r\right) \gamma $, $\widetilde{g}$: $[0,T]\times \mathbb{C\rightarrow
C}$.

We assume the following

${\mathbf{(H1)}}$. (Polynomial growth) \ $\widetilde{g}(s,\cdot )\in
C_{p}^{3}(\mathbb{C}\mathbf{,}\mathbb{C})$ and $h\in C_{p}^{3}(\mathbb{C}%
\mathbf{,}\mathbb{C}).$

${\mathbf{(H2)}}.$ (Bounded derivatives) $\alpha ,\beta ,\theta $ are $%
\mathbb{C}$-valued bounded and deterministic functions

${\mathbf{(H3)}}$. (Analyticity)\ $h,\widetilde{g}$ are analytic w.r.t
spatial variable. e.g. for $x\mathbb{=}x^{1}+ix^{2}$, $\frac{\partial h^{1}}{%
\partial x^{1}}=\frac{\partial h^{2}}{\partial x^{2}}$, $\frac{\partial h^{1}%
}{\partial x^{2}}=-\frac{\partial h^{2}}{\partial x^{1}}$.

\begin{remark}
From the Maximum Modulus Principle, for a nonlinear complex function $%
g\left( r,\mathbb{X},\mathbb{Y},\mathbb{Z},\mathbb{T}\right) $, the
derivatives of $g$ w.r.t $\mathbb{Y},\mathbb{Z},\mathbb{T}$ are bounded and
analytic on the whole complex plane means that $g$ is linear in $\mathbb{Y},%
\mathbb{Z},\mathbb{T}$.
\end{remark}

BSDE (3.1) is equivalent to the following 2-dimensional BSDE: for $s\in
\lbrack 0,T]$
\begin{eqnarray}
\left(
\begin{array}{l}
Y_{s}^{1} \\
Y_{s}^{2}%
\end{array}%
\right) &=&\left(
\begin{array}{l}
h^{1}X_{r}^{1},X_{r}^{2} \\
h^{2}\left( X_{r}^{1},X_{r}^{2}\right)%
\end{array}%
\right) +\int_{s}^{T}\left(
\begin{array}{l}
g^{1}\left(
r,X_{r}^{1},X_{r}^{2},Y_{r}^{1},Y_{r}^{2},Z_{r}^{1},Z_{r}^{2},\Gamma
_{r}^{1},\Gamma _{r}^{2}\right) \\
g^{2}\left(
r,X_{r}^{1},X_{r}^{2},Y_{r}^{1},Y_{r}^{2},Z_{r}^{1},Z_{r}^{2},\Gamma
_{r}^{1},\Gamma _{r}^{2}\right)%
\end{array}%
\right) dr  \notag \\
&&-\int_{s}^{T}\left(
\begin{array}{ll}
Z_{r}^{1} & -Z_{r}^{2} \\
Z_{r}^{2} & Z_{r}^{1}%
\end{array}%
\right) \left(
\begin{array}{l}
dB_{r}^{1} \\
dB_{r}^{2}%
\end{array}%
\right) -\int_{s}^{T}\left(
\begin{array}{ll}
\Gamma _{r}^{1} & \Gamma _{r}^{2} \\
\Gamma _{r}^{2} & -\Gamma _{r}^{1}%
\end{array}%
\right) \left(
\begin{array}{l}
dB_{r}^{1} \\
dB_{r}^{2}%
\end{array}%
\right) .
\end{eqnarray}%
or the following

\begin{equation}
\left(
\begin{array}{l}
Y_{s}^{1} \\
Y_{s}^{2}%
\end{array}%
\right) =\left(
\begin{array}{l}
h^{1}\left( X_{r}^{1},X_{r}^{2}\right) \\
h^{2}\left( X_{r}^{1},X_{r}^{2}\right)%
\end{array}%
\right) +\int_{s}^{T}\left(
\begin{array}{l}
f^{1}\left( r,\mathbf{X}_{r},\mathbf{Y}_{r},\mathbf{Z}_{r}\right) \\
f^{2}\left( r,\mathbf{X}_{r},\mathbf{Y}_{r},\mathbf{Z}_{r}\right)%
\end{array}%
\right) dr-\int_{s}^{T}\left(
\begin{array}{ll}
Z_{r}^{11} & Z_{r}^{21} \\
Z_{r}^{12} & Z_{r}^{22}%
\end{array}%
\right) \left(
\begin{array}{l}
dB_{r}^{1} \\
dB_{r}^{2}%
\end{array}%
\right) ,  \label{3.3}
\end{equation}%
where $\mathbf{X}=\left(
\begin{array}{l}
X^{1} \\
X^{2}%
\end{array}%
\right) $, $\mathbf{Y}_{s}=\left(
\begin{array}{l}
Y^{1} \\
Y^{2}%
\end{array}%
\right) $, $\mathbf{Z}=(\mathbf{Z}^{1},\mathbf{Z}^{2})=\left(
\begin{array}{ll}
Z^{11} & Z^{12} \\
Z^{21} & Z^{22}%
\end{array}%
\right) =\left(
\begin{array}{ll}
Z^{1}+\Gamma ^{1} & \Gamma ^{2}+Z^{2} \\
\Gamma ^{2}-Z^{2} & Z^{1}-\Gamma ^{1}%
\end{array}%
\right) $,

\begin{eqnarray}
f^{1}\left( r,\mathbf{X}_{r},\mathbf{Y}_{r},\mathbf{Z}_{r}\right) &=&g^{1}%
\mathbf{(}r,X_{r}^{1},X_{r}^{2},Y_{r}^{1},Y_{r}^{2},\frac{1}{2}%
(Z_{r}^{11}+Z_{r}^{22}),\frac{1}{2}(Z_{r}^{12}-Z_{r}^{21}),  \notag \\
&&\ \ \frac{1}{2}(Z_{r}^{11}-Z_{r}^{22}),\frac{1}{2}(Z_{r}^{12}+Z_{r}^{21})%
\mathbf{)},
\end{eqnarray}

\begin{eqnarray}
f^{2}\left( r,\mathbf{X}_{r},\mathbf{Y}_{r},\mathbf{Z}_{r}\right) &=&g^{2}%
\mathbf{(}r,X_{r}^{1},X_{r}^{2},Y_{r}^{1},Y_{r}^{2},\frac{1}{2}%
(Z_{r}^{11}+Z_{r}^{22}),\frac{1}{2}(Z_{r}^{12}-Z_{r}^{21}),  \notag \\
&&\ \ \frac{1}{2}(Z_{r}^{11}-Z_{r}^{22}),\frac{1}{2}(Z_{r}^{12}+Z_{r}^{21})%
\mathbf{)}.
\end{eqnarray}

One can check that $f=(f^{1},f^{2})$ is analytic in $\left( \mathbf{X},%
\mathbf{Y},\mathbf{Z}\right) $, e.g. $\frac{\partial f^{1}}{\partial \mathbf{%
Z}^{1}}=\frac{\partial f^{2}}{\partial \mathbf{Z}^{2}}$, $\frac{\partial
f^{1}}{\partial \mathbf{Z}^{2}}=-\frac{\partial f^{2}}{\partial \mathbf{Z}%
^{1}}$, and satisfies the Lipschitz condition in $\left( \mathbf{Y},\mathbf{Z%
}\right) $, thus there is a unique pair $\left( \mathbf{Y}_{s}^{t,x},\mathbf{%
Z}_{s}^{t,x}\right) \in {\mathcal{S}}_{\mathcal{F}}^{\mathrm{2}}\left( t,{%
\mathit{T}};\mathbf{R}^{2}\right) \times {\mathcal{H}_{\mathcal{F}}^{\mathrm{%
2}}}(t,T;\mathbf{R}^{2}\times \mathbf{R}^{2})$ which solves the real BSDE
(3.3) (see Pardoux and Peng 1990). Therefore there is a unique triple $%
\left( \mathbb{Y}_{s}^{t,x},\mathbb{Z}_{s}^{t,x},\mathbb{T}_{s}^{t,x}\right)
\in {\mathcal{S}}_{\mathcal{F}}^{\mathrm{2}}\left( t,{\mathit{T}};\mathbb{C}%
\right) \times {\mathcal{H}_{\mathcal{F}}^{\mathrm{2}}}(t,T;\mathbb{C}%
)\times {\mathcal{H}_{\mathcal{F}}^{\mathrm{2}}}(t,T;\mathbb{C})$ for BSDE
(3.1).

\begin{remark}
The analyticity is not used for the existence and uniqueness of the
solutions for BSDE (3.1). It is just useful when we derive the Feynman-Kac
formula.
\end{remark}

\begin{theorem}
$\{\mathbb{Y}_s^{t,x};t\leq s\leq T,x\in \mathbb{C}\}$ is analytic in $x$
and continuous in $(s,t)$.
\end{theorem}

Before proceeding to the proof, we first state a useful corollary:

\begin{corollary}
For any $t\in [0,T]$, the mapping $x\rightarrow \mathbb{Y}_t^{t,x}$ is
analytic, the function and its partial derivatives of order one and two
being continuous in $(t,x)$.
\end{corollary}

\textbf{Proof of Theorem 4.1.} It suffice to prove the following

Step 1. $\{\mathbb{Y}_s^{t,x}\}$ is continuous in $(s,t)$, $\{\frac{\partial
\mathbb{Y}_s^{t,x}}{\partial x}\}$ is continuous in $x$.

Step 2. Cauchy-Riemann equations: $\frac{\partial Y_s^1}{\partial x^1}=\frac{%
\partial Y_s^2}{\partial x^2}$, $\frac{\partial Y_s^1}{\partial x^2}=-\frac{%
\partial Y_s^2}{\partial x^1}$.

Step 1 follows immediately the fact that

\begin{equation*}
{\mathbf{E}}\left[ {\mathrm{sup}}_{t\wedge \hat{t}\leq s\leq T}|\mathbf{Y}%
_{s}^{t,x}-\mathbf{Y}_{s}^{\hat{t},\hat{x}}|^{p}\right] \leq
c_{p}(1+|x|^{q})(|x-\hat{x}|^{p}+|t-\hat{t}|^{\frac{p}{2}})
\end{equation*}%
\textbf{\ }
\begin{equation*}
{\mathbf{E}}\left[ {\mathrm{sup}}_{t\wedge \hat{t}\leq s\leq
T}|\bigtriangleup _{h}^{i}\mathbf{Y}_{s}^{t,x}-\bigtriangleup _{h}^{i}%
\mathbf{Y}_{s}^{\hat{t},\hat{x}}|^{p}\right] \leq c_{p}(1+|x|^{q}+|\hat{x}%
|^{q}+|h|^{q}+|\hat{h}|^{q})(|x-\hat{x}|^{p}+|h-\hat{h}|^{p}+|t-\hat{t}|^{%
\frac{p}{2}})
\end{equation*}%
where $\bigtriangleup _{h}^{i}\mathbf{Y}_{s}^{t,x}=(\mathbf{Y}%
_{s}^{t,x+he^{i}}-\mathbf{Y}_{s}^{t,x})/h$, $h\in \mathbf{R\backslash \{}0%
\mathbf{\}}$\textbf{, }$\mathbf{\{}e^{1},e^{2}\mathbf{\}}$ is an orthogonal
basis of $\mathbf{R}^{2}$.

We now prove the Cauchy-Riemann equations. Pardoux and Peng (1992, eq. 13)
says that,
\begin{eqnarray}
\nabla \mathbf{Y}_s^{t,x} &=&\int_s^T[f_x^{\prime }\left( r,\mathbf{X}%
_r^{t,x},\mathbf{Y}_r^{t,x},\mathbf{Z}_r^{t,x}\right) \nabla \mathbf{X}%
_r^{t,x}+f_y^{\prime }\left( r,\mathbf{X}_r^{t,x},\mathbf{Y}_r^{t,x},\mathbf{%
Z}_r^{t,x}\right) \nabla \mathbf{Y}_r^{t,x}  \notag \\
&&+f_z^{\prime }\left( r,\mathbf{X}_r^{t,x},\mathbf{Y}_r^{t,x},\mathbf{Z}%
_r^{t,x}\right) \nabla \mathbf{Z}_r^{t,x}]dr+\int_s^T\left\langle (\nabla
\mathbf{Z}_r^{t,x})^{*},d\mathbf{B}_r\right\rangle ,s\in [t,T].
\end{eqnarray}
where $\nabla \mathbf{Y}_s^{t,x}$ is the matrix of first order partial
derivatives of $\mathbf{Y}_s^{t,x}$ ($x$ denotes the initial condition of
SDE (2.1)). $\nabla \mathbf{X}_s^{t,x}$ and $\nabla \mathbf{Z}_s^{t,x}$ are
defined analogously. Let $(\mathbf{M}_t^s)_{t\leq s\leq T}$ be the solution
of the following matrix-valued SDE:

\begin{eqnarray}
d\mathbf{M}_s &=&\mathbf{M}_sf_y^{\prime }ds+\left\langle \mathbf{M}%
_sf_z^{\prime },d\mathbf{B}_s\right\rangle ,t\leq s\leq T.  \notag \\
\mathbf{M}_t &=&\mathbf{I}
\end{eqnarray}
which is given by
\begin{equation}
\mathbf{M}_t^s=\exp \{\int_t^s\left[ f_y^{\prime }-\frac 12f_z^{\prime
}\cdot (f_z^{\prime })^{*}\right] dr+\int_t^s\left\langle f_z^{\prime },d%
\mathbf{B}_r\right\rangle \},t\leq s\leq T.
\end{equation}
Since $f_y^{\prime }$, $f_z^{\prime }$ is of class $\mathbb{C}_L$, $\mathbf{M%
}_t^s$ is also of class $\mathbb{C}_L$.

Then applying a 4-dimensional version of It\^{o}'s formula (see \textsc{$%
\emptyset $ksendal} 2005, Th.4.2.1) to $\mathbf{M}_{s}\nabla \mathbf{Y}%
_{s}^{t,x}$, we deduce that

\begin{eqnarray*}
d(\mathbf{M}_{s}\nabla \mathbf{Y}_{s}^{t,x}) &=&\mathbf{M}_{s}d(\nabla
\mathbf{Y}_{s}^{t,x})+(d\mathbf{M}_{s})\nabla \mathbf{Y}_{s}^{t,x}+(d\mathbf{%
M}_{s})_{s}d(\nabla \mathbf{Y}_{s}^{t,x}) \\
&=&-\mathbf{M}_{s}\left( f_{y}^{\prime }\nabla \mathbf{Y}_{s}^{t,x}+f_{z}^{%
\prime }\nabla \mathbf{Z}_{s}^{t,x}+f_{x}^{\prime }\nabla \mathbf{X}%
_{s}^{t,x}\right) ds+\left\langle \mathbf{M}_{s}\nabla \mathbf{Z}_{s}^{t,x},d%
\mathbf{B}_{s}\right\rangle \\
&&+\mathbf{M}_{s}f_{y}^{\prime }\nabla \mathbf{Y}_{s}^{t,x}ds+\left\langle
\mathbf{M}_{s}f_{z}^{\prime }\nabla \mathbf{Y}_{s}^{t,x},d\mathbf{B}%
_{s}\right\rangle +\mathbf{M}_{s}f_{z}^{\prime }\nabla \mathbf{Z}%
_{s}^{t,x}ds,s\in \lbrack t,T].
\end{eqnarray*}%
Therefore
\begin{equation*}
\mathbf{M}_{s}\nabla \mathbf{Y}_{s}^{t,x}=\mathbf{M}_{T}\nabla \mathbf{Y}%
_{T}^{t,x}+\int_{s}^{T}\mathbf{M}_{r}f_{x}^{\prime }\nabla \mathbf{X}%
_{r}^{t,x}dr-\int_{s}^{T}\left\langle \mathbf{M}_{r}(f_{z}^{\prime }\nabla
\mathbf{Y}_{r}^{t,x}+\nabla \mathbf{Z}_{r}^{t,x})^{\ast },d\mathbf{B}%
_{r}\right\rangle ,
\end{equation*}%
then

\begin{equation}
\nabla \mathbf{Y}_{s}^{t,x}={\normalsize E}\left[ \mathbf{M}_{s}^{T}{\large h%
}^{\prime }\left( \mathbf{X}_{T}^{t,x}\right) \nabla \mathbf{X}%
_{T}^{t,x}+\int_{s}^{T}\mathbf{M}_{s}^{r}f_{x}^{\prime }\nabla \mathbf{X}%
_{r}^{t,x}dr\mid \mathcal{F}_{r}\right]
\end{equation}%
By Lemma \ref{lem2.2}, It is known that $\nabla \mathbf{Y}_{s}^{t,x}$ is of
class $\mathbb{C}_{L}$, thus the Cauchy-Riemann equations hold true. The
proof is complete. $\Box $

\begin{remark}
Let $\mathbf{Z}_s^{t,x}=(Z_s^1,Z_s^2)=\left(
\begin{array}{ll}
Z_s^{11} & Z_s^{12} \\
Z_s^{21} & Z_s^{22}%
\end{array}
\right) $, then $f_z^{\prime }=(A,B)$, $(f_z^{\prime })^{*}=\left(
\begin{array}{l}
A \\
B%
\end{array}
\right) $, $\nabla \mathbf{Z}_s^{t,x}=\left(
\begin{array}{l}
C \\
D%
\end{array}
\right) $, $(\nabla \mathbf{Z}_s^{t,x})^{*}=(C,D)$, where $A=\left(
\begin{array}{ll}
\frac{\partial f^1}{\partial Z^{11}} & \frac{\partial f^1}{\partial Z^{12}}
\\
\frac{\partial f^2}{\partial Z^{11}} & \frac{\partial f^2}{\partial Z^{12}}%
\end{array}
\right) $, $B=\left(
\begin{array}{ll}
\frac{\partial f^1}{\partial Z^{21}} & \frac{\partial f^1}{\partial Z^{22}}
\\
\frac{\partial f^2}{\partial Z^{21}} & \frac{\partial f^2}{\partial Z^{22}}%
\end{array}
\right) $, $C=\left(
\begin{array}{ll}
\frac{\partial Z^{11}}{\partial x^1} & \frac{\partial Z^{11}}{\partial x^2}
\\
\frac{\partial Z^{12}}{\partial x^1} & \frac{\partial Z^{11}}{\partial x^2}%
\end{array}
\right) $, $D=\left(
\begin{array}{ll}
\frac{\partial Z^{21}}{\partial x^1} & \frac{\partial Z^{21}}{\partial x^2}
\\
\frac{\partial Z^{22}}{\partial x^1} & \frac{\partial Z^{22}}{\partial x^2}%
\end{array}
\right) .$

Since $f$ is analytic in $\mathbf{Z}=Z^1+i\cdot Z^2$, one can check that $%
A,B $ are of class $\mathbb{C}_L$.
\end{remark}

For BSDE \eqref{3.3}, Pardoux and Peng (1992) proved that: for any $\ 0\leq
t\leq s\leq T$, $x\in \mathbf{R}^{n}$, $(\mathbf{Z}_{s}^{t,x})^{\ast
}=\nabla \mathbf{Y}_{s}^{t,x}(\nabla \mathbf{X}_{s}^{t,x})^{-1}\sigma
_{s}^{\ast }$ and particularly $(\mathbf{Z}_{t}^{t,x})^{\ast }=\nabla
\mathbf{Y}_{t}^{t,x}\sigma _{t}^{\ast }$. Note that $\nabla \mathbf{X}%
_{s}^{t,x}=\mathbf{I}$ in this paper, for BSDE \eqref{3.1}, we have the
following results.

\begin{proposition}
\label{prop3.1}For$\ 0\leq t\leq s\leq T$, $x\in \mathbb{C}$,

\begin{equation}
(\mathbb{Z}_{s}^{t,x})=\frac{d\mathbb{Y}_{s}^{t,x}}{dx}\sigma _{s},
\end{equation}%
\begin{equation}
(\mathbb{T}_{s}^{t,x})=\frac{d\mathbb{Y}_{s}^{t,x}}{dx}\gamma _{s},
\end{equation}
\end{proposition}

\textbf{Proof}. By Pardoux and Peng (1992, Lemma 2.5), we get that

\begin{equation}
(\mathbf{Z}^{t,x})^{\ast }=\left(
\begin{array}{ll}
Z^{1}+\Gamma ^{1} & \Gamma ^{2}-Z^{2} \\
\Gamma ^{2}+Z^{2} & Z^{1}-\Gamma ^{1}%
\end{array}%
\right) =\left(
\begin{array}{ll}
\frac{\partial Y^{1}}{\partial x^{1}} & \frac{\partial Y^{1}}{\partial x^{2}}
\\
\frac{\partial Y^{2}}{\partial x^{1}} & \frac{\partial Y^{2}}{\partial x^{2}}%
\end{array}%
\right) \left(
\begin{array}{ll}
\sigma ^{1}+\gamma ^{1} & \gamma ^{2}-\sigma ^{2} \\
\gamma ^{2}+\sigma ^{2} & \sigma ^{1}-\gamma ^{1}%
\end{array}%
\right) .
\end{equation}%
From the above equation, we have that

\begin{equation}
\left\{ \begin{aligned} Z^1+\Gamma ^1 &= \frac{\partial Y^1}{\partial
x^1}(\sigma ^1+\gamma ^1)+ \frac{\partial Y^1}{\partial x^2}(\sigma ^2+
\gamma^2), \\ Z^1-\Gamma ^1&=\frac{\partial Y^1}{\partial x^1}(\sigma
^1-\gamma ^1)+\frac{ \partial Y^1}{\partial x^2}(\sigma ^2- \gamma^2),
\end{aligned} \right.
\end{equation}
and

\begin{equation}
\left\{ \begin{aligned} \Gamma ^2-Z^2 &= \frac{\partial Y^2}{\partial
x^1}(\gamma ^1-\sigma ^1)+ \frac{\partial Y^2}{\partial x^2}(\gamma
^2-\sigma ^2), \\ \Gamma ^2+Z^2&=\frac{\partial Y^2}{\partial x^1}(\gamma
^1+\sigma ^1)+\frac{ \partial Y^2}{\partial x^2}(\gamma ^2+\sigma ^2).
\end{aligned} \right.
\end{equation}
Therefore,
\begin{equation}
\left\{ \begin{aligned} Z^1&= \frac{\partial Y^1}{\partial x^1}\sigma ^1+
\frac{\partial Y^1}{\partial x^2}\sigma ^2, \\ \Gamma ^1&=\frac{\partial
Y^1}{\partial x^1}\gamma ^1+\frac{ \partial Y^1}{\partial x^2}\gamma^2,
\end{aligned} \right.
\end{equation}
and

\begin{equation}
\left\{ \begin{aligned} Z^2 &= \frac{\partial Y^2}{\partial x^1}\sigma ^1+
\frac{\partial Y^2}{\partial x^2}\sigma ^2, \\ \Gamma ^2&=\frac{\partial
Y^2}{\partial x^1}\gamma ^1+\frac{ \partial Y^2}{\partial x^2}\gamma ^2.
\end{aligned}\right.
\end{equation}%
Thus
\begin{equation*}
Z^{1}+i\cdot Z^{2}=(\frac{\partial Y^{1}}{\partial x^{1}}\sigma ^{1}+\frac{%
\partial Y^{1}}{\partial x^{2}}\sigma ^{2})+i(\frac{\partial Y^{2}}{\partial
x^{1}}\sigma ^{1}+\frac{\partial Y^{2}}{\partial x^{2}}\sigma ^{2})=(\frac{%
\partial Y^{1}}{\partial x^{1}}+i\frac{\partial Y^{2}}{\partial x^{1}}%
)(\sigma ^{1}+i\sigma ^{2}),
\end{equation*}

\begin{equation*}
\Gamma ^{1}+i\cdot \Gamma ^{2}=(\frac{\partial Y^{1}}{\partial x^{1}}\gamma
^{1}+\frac{\partial Y^{1}}{\partial x^{2}}\gamma ^{2})+i(\frac{\partial Y^{2}%
}{\partial x^{1}}\gamma ^{1}+\frac{\partial Y^{2}}{\partial x^{2}}\gamma
^{2})=(\frac{\partial Y^{1}}{\partial x^{1}}+i\frac{\partial Y^{2}}{\partial
x^{1}})(\gamma ^{1}+i\gamma ^{2}),
\end{equation*}%
that is

\begin{eqnarray*}
(\mathbb{Z}_{s}^{t,x}) &=&\frac{d\mathbb{Y}_{t}^{t,x}}{dx}\ \sigma _{s}, \\
(\mathbb{T}_{s}^{t,x}) &=&\frac{d\mathbb{Y}_{t}^{t,x}}{dx}\ \gamma _{s}.
\end{eqnarray*}

$\Box $

By Proposition \ref{prop3.1} we know that, if $\sigma =0$(resp. $\gamma =0$%
), then $\mathbb{Z}=0$(resp. $\mathbb{T}=0$). Since there is a unique
solution ($\mathbb{Y},\mathbb{Z},\mathbb{T}$) for BSDE (3.1), we have the
following results:

\begin{corollary}
Under conditions (\textbf{H1})\textit{\ }$\thicksim $(\textbf{H3}), there is
a unique solution in ${\mathcal{S}}_{\mathcal{F}}^{\mathrm{2}}\left( t,{%
\mathit{T}};\mathbf{R}^{2}\right) \times {\mathcal{H}_{\mathcal{F}}^{\mathrm{%
2}}}(t,T;\mathbf{R}^{2}\times \mathbf{R}^{2})$ respectively for the
following two real forward-backward SDEs:%
\begin{equation}
\left\{
\begin{array}{ll}
\left(
\begin{array}{l}
X_{s}^{1} \\
X_{s}^{2}%
\end{array}%
\right) =\left(
\begin{array}{l}
x^{1} \\
x^{2}%
\end{array}%
\right) +\int_{t}^{s}\left(
\begin{array}{ll}
\sigma _{r}^{1} & -\sigma _{r}^{2} \\
\sigma _{r}^{2} & \sigma _{r}^{1}%
\end{array}%
\right) \left(
\begin{array}{l}
dB_{r}^{1} \\
dB_{r}^{2}%
\end{array}%
\right) \text{,} & s\in \left[ t,T\right] \text{,} \\
\begin{array}{l}
\left(
\begin{array}{l}
Y_{s}^{1} \\
Y_{s}^{2}%
\end{array}%
\right) =\left(
\begin{array}{l}
h^{1}\left( X_{r}^{1},X_{r}^{2}\right) \\
h^{2}\left( X_{r}^{1},X_{r}^{2}\right)%
\end{array}%
\right) +\int_{s}^{T}\left(
\begin{array}{l}
g^{1}\left(
r,X_{r}^{1},X_{r}^{2},Y_{r}^{1},Y_{r}^{2},Z_{r}^{1},Z_{r}^{2},\Gamma
_{r}^{1},\Gamma _{r}^{2}\right) \\
g^{2}\left(
r,X_{r}^{1},X_{r}^{2},Y_{r}^{1},Y_{r}^{2},Z_{r}^{1},Z_{r}^{2},\Gamma
_{r}^{1},\Gamma _{r}^{2}\right)%
\end{array}%
\right) dr \\
\ \ \ \ \ \ \ \ \ \ \ \ \ \ -\int_{s}^{T}\left(
\begin{array}{ll}
Z_{r}^{1} & -Z_{r}^{2} \\
Z_{r}^{2} & Z_{r}^{1}%
\end{array}%
\right) \left(
\begin{array}{l}
dB_{r}^{1} \\
dB_{r}^{2}%
\end{array}%
\right) \text{,}%
\end{array}
& s\in \left[ t,T\right] \text{,}%
\end{array}%
\right.
\end{equation}%
and%
\begin{equation}
\left\{
\begin{array}{ll}
\left(
\begin{array}{l}
X_{s}^{1} \\
X_{s}^{2}%
\end{array}%
\right) =\left(
\begin{array}{l}
x^{1} \\
x^{2}%
\end{array}%
\right) +\int_{t}^{s}\left(
\begin{array}{ll}
\gamma _{r}^{1} & \gamma _{r}^{2} \\
\gamma _{r}^{2} & -\gamma _{r}^{1}%
\end{array}%
\right) \left(
\begin{array}{l}
dB_{r}^{1} \\
dB_{r}^{2}%
\end{array}%
\right) \text{,} & s\in \left[ t,T\right] \text{,} \\
\begin{array}{l}
\left(
\begin{array}{l}
Y_{s}^{1} \\
Y_{s}^{2}%
\end{array}%
\right) =\left(
\begin{array}{l}
h^{1}\left( X_{r}^{1},X_{r}^{2}\right) \\
h^{2}\left( X_{r}^{1},X_{r}^{2}\right)%
\end{array}%
\right) +\int_{s}^{T}\left(
\begin{array}{l}
g^{1}\left(
r,X_{r}^{1},X_{r}^{2},Y_{r}^{1},Y_{r}^{2},Z_{r}^{1},Z_{r}^{2},\Gamma
_{r}^{1},\Gamma _{r}^{2}\right) \\
g^{2}\left(
r,X_{r}^{1},X_{r}^{2},Y_{r}^{1},Y_{r}^{2},Z_{r}^{1},Z_{r}^{2},\Gamma
_{r}^{1},\Gamma _{r}^{2}\right)%
\end{array}%
\right) dr \\
\ \ \ \ \ \ \ \ \ \ \ \ \ \ -\int_{s}^{T}\left(
\begin{array}{ll}
\Gamma _{r}^{1} & \Gamma _{r}^{2} \\
\Gamma _{r}^{2} & -\Gamma _{r}^{1}%
\end{array}%
\right) \left(
\begin{array}{l}
dB_{r}^{1} \\
dB_{r}^{2}%
\end{array}%
\right) \text{,}%
\end{array}
& s\in \left[ t,T\right] \text{,}%
\end{array}%
\right.
\end{equation}
\end{corollary}

\begin{remark}
The above two BSDEs are real-valued BSDEs with $Z$-constraints. Usually
there are no solutions for constrained BSDEs.
\end{remark}

\section{Complex BSDE and associated PDE}

Consider the following complex PDE:

\begin{equation}
\left\{ \begin{aligned} \Bbb{U}_t(t,x) &= -g\left( t,x,\Bbb{U}(t,x),\sigma
_t \Bbb{U}_x(t,x),\gamma _t \Bbb{U}_x(t,x)\right), \\ \Bbb{U}(T,x)&=h(x),\
0\leq t\leq T. \end{aligned}\right.  \label{4.1}
\end{equation}%
where $\mathbb{U}:\mathbf{R}^{+}\times \mathbb{C}\rightarrow \mathbb{C}$, $%
\sigma =\sigma ^{1}+i\sigma ^{2}$, $\gamma =\gamma ^{1}+i\gamma ^{2}$, $%
\sigma ^{1}\gamma ^{1}=\sigma ^{2}\gamma ^{2}$.

\begin{theorem}
Let $h,g$ satisfy (\textbf{H1})\textit{\ }$\thicksim $(\textbf{H3}). If for
each (t,x), $\{\mathbb{U}(t,x);t\leq s\leq T,x\in \mathbb{C}\}$ is analytic
with respect to $x$ and continuous in $t$, and satisfies PDE (4.1), then
\begin{equation*}
(\mathbb{Y}_{s}^{t,x},\mathbb{Z}_{s}^{t,x},\mathbb{T}_{s}^{t,x}):=\left(
\mathbb{U}(s,\mathbb{X}_{s}^{t,x}),\sigma _{s}\mathbb{U}_{x}(s,\mathbb{X}%
_{s}^{t,x}),\gamma _{s}\mathbb{U}_{x}(s,\mathbb{X}_{s}^{t,x})\right)
\end{equation*}%
solves BSDE (3.1). Furthermore,
\begin{equation*}
\left( \mathbb{U}(t,x),\sigma _{t}\mathbb{U}_{x}(t,x),\gamma _{t}\mathbb{U}%
_{x}(t,x)\right) _{0\leq t\leq T}=\left( \mathbb{Y}_{t}^{t,x},\mathbb{Z}%
_{t}^{t,x},\mathbb{T}_{t}^{t,x}\right) ,
\end{equation*}%
and
\begin{equation*}
\mathbb{U}(t,x)=\mathbb{Y}_{t}^{t,x}=E\left[ h\left( \mathbb{X}%
_{T}^{t,x}\right) +\int_{t}^{T}g\left( r,\mathbb{X}_{r}^{t,x},\mathbb{Y}%
_{r}^{t,x},\mathbb{Z}_{r}^{t,x},\mathbb{T}_{r}^{t,x}\right) dr\right] .
\end{equation*}
\end{theorem}

Before proving the above theorem, we need an It\^{o}'s lemma.

\begin{lemma}
Let $\mathbb{X}_{t}\in {\mathcal{S}}_{\mathcal{F}}^{\mathrm{2}}\left( 0,{%
\mathit{T}};\mathbb{C}\right) $, $b_{t},\sigma _{t},\gamma _{t}\times {%
\mathcal{H}_{\mathcal{F}}^{\mathrm{2}}}(0,{\mathit{T}};\mathbb{C})$, such
that $\sigma ^{1}\gamma ^{1}=\sigma ^{2}\gamma ^{2}$ for almost all $t$ and
\begin{equation}
d\mathbb{X}_{t}=b_{t}dt+\sigma _{t}d\mathbb{B}_{t}+\gamma _{t}d\bar{\mathbb{B%
}}_{t},t\geq 0.
\end{equation}

If $\mathbb{F}(t,x)=u(t,x^{1},x^{2})+iv(t,x^{1},x^{2})$ is an analytic
function w.r.t the complex variable x and continuous in t, i.e. for any t, $%
\forall x\in \mathbb{C}$, $\mathbb{F}$ satisfies the Cauchy-Riemann
equations: $\frac{\partial u}{\partial x^{1}}=\frac{\partial v}{\partial
x^{2}}$, $\frac{\partial u}{\partial x^{2}}=-\frac{\partial v}{\partial x^{1}%
}$, then
\begin{eqnarray}
d\mathbb{F}(t,\mathbb{X}_{t}) &=&\frac{\partial \mathbb{F}}{\partial t}(t,%
\mathbb{X}_{t})dt+\frac{\partial \mathbb{F}}{\partial x}(t,\mathbb{X}_{t})d%
\mathbb{X}_{t}  \notag \\
&=&\frac{\partial \mathbb{F}}{\partial t}(t,\mathbb{X}_{t})dt+\frac{\partial
\mathbb{F}}{\partial x}(t,\mathbb{X}_{t})b_{t}dt+\frac{\partial \mathbb{F}}{%
\partial x}(t,\mathbb{X}_{t})\sigma _{t}d\mathbb{B}_{t}+\frac{\partial
\mathbb{F}}{\partial x}(t,\mathbb{X}_{t})\gamma _{t}d\bar{\mathbb{B}}%
_{t},t\geq 0.
\end{eqnarray}%
where $\frac{\partial \mathbb{F}}{\partial x}(t,x)$ is the complex partial
derivative of $\mathbb{F}$ w.r.t $x$.
\end{lemma}

Proof.
\begin{eqnarray*}
d\mathbb{F}(t,\mathbb{X}_{t})
&=&du(t,X_{t}^{1},X_{t}^{2})+idv(t,X_{t}^{1},X_{t}^{2}) \\
&=&\frac{\partial u}{\partial x^{1}}dX_{t}^{1}+\frac{\partial u}{\partial
x^{2}}dX_{t}^{2}+\frac{\partial u}{\partial t}dt+\frac{1}{2}(\frac{\partial
^{2}u}{\partial (x^{1})^{2}}d<X^{1}>_{t}+\frac{\partial ^{2}u}{\partial
(x^{2})^{2}}d<X^{2}>_{t}) \\
&&+i\left[ \frac{\partial v}{\partial x^{1}}dX_{t}^{1}+\frac{\partial v}{%
\partial x^{2}}dX_{t}^{2}+\frac{\partial v}{\partial t}dt+\frac{1}{2}(\frac{%
\partial ^{2}v}{\partial (x^{1})^{2}}d<X^{1}>_{t}+\frac{\partial ^{2}v}{%
\partial (x^{2})^{2}}d<X^{2}>_{t})\right] \\
&=&\left( \frac{\partial u}{\partial x^{1}}+i\frac{\partial v}{\partial x^{1}%
}\right) (dX_{t}^{1}+idX_{t}^{2})+\left( \frac{\partial u}{\partial t}+i%
\frac{\partial v}{\partial t}\right) dt \\
&=&\frac{\partial \mathbb{F}}{\partial t}(t,\mathbb{X}_{t})dt+\frac{\partial
\mathbb{F}}{\partial x}(t,\mathbb{X}_{t})d\mathbb{X}_{t}
\end{eqnarray*}

where we have used the conjugate harmonicity of function $u$, $v$ and the
condition $\sigma ^{1}\gamma ^{1}=\sigma ^{2}\gamma ^{2}$. $\Box $

\textbf{Proof of Theorem 4.1}. It suffices to show that
\begin{equation*}
\{\mathbb{U}(s,\mathbb{X}_{s}^{t,x}),\sigma _{s}\mathbb{U}_{x}(s,\mathbb{X}%
_{s}^{t,x}),\gamma _{s}\mathbb{U}_{x}(s,\mathbb{X}_{s}^{t,x});t\leq s\leq T\}
\end{equation*}%
solves BSDE (3.1). Applying the complex It\^{o} formula to $\mathbb{U}(s,%
\mathbb{X}_{s}^{t,x})$ between $s=t$ and $s=T$, we get that
\begin{eqnarray*}
\mathbb{Y}_{s}^{t,x} &=&h\left( \mathbb{X}_{T}^{t,x}\right) +\int_{s}^{T}%
\left[ g\left( r,\mathbb{X}_{r}^{t,x},\mathbb{U}(r,\mathbb{X}%
_{r}^{t,x}),\sigma _{r}\mathbb{U}_{x}(r,\mathbb{X}_{r}^{t,x}),\gamma _{r}%
\mathbb{U}_{x}(r,\mathbb{X}_{r}^{t,x})\right) \right] dr \\
&&-\int_{s}^{T}\sigma _{r}\mathbb{U}_{x}(r,\mathbb{X}_{r}^{t,x})d\mathbb{B}%
_{r}-\int_{s}^{T}{}\gamma _{r}\mathbb{U}_{x}(r,\mathbb{X}_{r}^{t,x})d\bar{%
\mathbb{B}}_{r},\ t\leq s\leq T.
\end{eqnarray*}

Thus $(\mathbb{Y}_{s}^{t,x},\mathbb{Z}_{s}^{t,x},\mathbb{T}_{s}^{t,x})=\{%
\mathbb{U}(s,\mathbb{X}_{s}^{t,x}),\sigma _{s}\mathbb{U}_{x}(s,\mathbb{X}%
_{s}^{t,x}),\gamma _{s}\mathbb{U}_{x}(s,\mathbb{X}_{s}^{t,x})\}$ solves BSDE
(3.1). $\Box $

We now show the converse of Theorem 4.1.

\begin{theorem}
\label{thm4.3}Let $h,g$ satisfy (\textbf{H1})\textit{\ }$\thicksim $(\textbf{%
H3}). Let $\{(\mathbb{Y}_{s}^{t,x});t\leq s\leq T\}$ be the solution of BSDE
(3.1). Then $\left( \mathbb{U}(t,x)\right) _{0\leq t\leq T}=\left( \mathbb{Y}%
_{t}^{t,x}\right) _{0\leq t\leq T}$ is the unique classical solution of
backward PDE (4.1) and $\mathbb{U}(t,x)$ is analytic in $x$ for each $t$.
\end{theorem}

\textbf{Proof}. Uniqueness follows from Theorem 4.1. We now prove that $%
\left( \mathbb{Y}_{t}^{t,x}\right) $ is a solution to PDE (4.1). Let $\delta
>0$ s.t. $t+\delta \leq T$. Clearly $\mathbb{Y}_{t+\delta }^{t,x}=\mathbb{Y}%
_{t+\delta }^{t+\delta ,\mathbb{X}_{t+\delta }^{t,x}}$. Hence by the complex
It\^{o} formula and the analyticity of $\mathbb{U}(t,x)$ in $x$ and BSDE %
\eqref{3.1}, we have

\begin{eqnarray*}
\mathbb{U}(t+\delta ,x)-\mathbb{U}(t,x) &=&\left[ \mathbb{U}(t+\delta ,x)-%
\mathbb{U}(t+\delta ,\mathbb{X}_{t+\delta }^{t,x})\right] +\left[ \mathbb{U}%
(t+\delta ,\mathbb{X}_{t+\delta }^{t,x})-\mathbb{U}(t,x)\right] \\
&=&-\int_{t}^{t+\delta }g\left( r,x,\mathbb{X}_{r}^{t,x},\mathbb{Y}%
_{r}^{t,x},\mathbb{Z}_{r}^{t,x},\mathbb{T}_{r}^{t,x}\right) dr \\
&&-\int_{t}^{t+\delta }\sigma _{r}\mathbb{U}_{x}(t+\delta ,\mathbb{X}%
_{r}^{t,x})d\mathbb{B}_{r} \\
&&-\int_{t}^{t+\delta }\gamma _{r}\mathbb{U}_{x}(t+\delta ,\mathbb{X}%
_{r}^{t,x})d\bar{\mathbb{B}}_{r} \\
&&+\int_{t}^{t+\delta }\mathbb{Z}_{r}^{t,x}d\mathbb{B}_{r}+\int_{t}^{t+%
\delta }\mathbb{T}_{r}^{t,x}d\bar{\mathbb{B}}_{r},
\end{eqnarray*}

Let $t=t_{0}\leq t_{1}\leq ...\leq t_{n}=T$, we get
\begin{eqnarray}
h(x)-\mathbb{U}(t,x) &=&-\sum_{i=0}^{n-1}\int_{t_{i}}^{t_{i+1}}g\left( r,x,%
\mathbb{X}_{r}^{t,x},\mathbb{Y}_{r}^{t,x},\mathbb{Z}_{r}^{t,x},\mathbb{T}%
_{r}^{t,x}\right) dr  \notag \\
&&+\int_{t_{i}}^{t_{i+1}}\left[ \mathbb{Z}_{r}^{t,x}-\sigma _{r}\mathbb{U}%
_{x}(t_{i+1},\mathbb{X}_{r}^{t,x})\right] d\mathbb{B}_{r}  \notag \\
&&+\int_{t_{i}}^{t_{i+1}}\left[ \mathbb{T}_{r}^{t,x}-\gamma _{r}\mathbb{U}%
_{x}(t_{i+1},\mathbb{X}_{r}^{t,x})\right] d\bar{\mathbb{B}}_{r}.
\end{eqnarray}

Let the mesh size ${\mathrm{sup}}_{0\leq i\leq n-1}$ $(t_{i+1}-t_i)%
\rightarrow 0$, we obtain the limit

\begin{equation*}
\mathbb{U}(t,x)=h(x)+\int_{s}^{T}g\left( r,x,\mathbb{U}(r,x),\sigma _{r}%
\mathbb{U}_{x}(r,x),\gamma _{r}\mathbb{U}_{x}(r,x)\right) dr,\ t\leq s\leq T.
\end{equation*}%
$\Box $

\begin{remark}
\label{app} Viscosity solution for PDE \eqref{4.1} is not involved in the
present paper because, the analyticity of parameters of PDE (4.1) leads to
existence of the first and the second order derivatives.
\end{remark}

The well-known Cauchy--Kovalevski theorem states a local existence and
uniqueness of solution for partial differential equations whose coefficients
are analytic functions, associated with Cauchy initial value problems. A
special case was proven by Cauchy in 1842, and the full result by Kowalevski
(1875). Theorem \ref{thm4.3} extends the first order Cauchy--Kovalevski
theorem to the case of global solutions.

\end{document}